\documentclass[12pt]{article} 
\usepackage{latexsym,amsmath}
\usepackage{amssymb,graphicx} 

\newtheorem{theo}{Theorem}[section] 
\newtheorem{pro}{Proposition}[section] 
\newtheorem{lem}{Lemma}[section]

\newcommand{\qedbox}{ \fbox{}}
\newenvironment{proof}{\noindent\textsc{Proof: }}{\hfill$\qedbox$}
\newenvironment{proofmaintheo}{\noindent\textsc{Proof of the Theorem \ref{maintheo}: }}{\hfill$\qedbox$}

\newenvironment{proofnormi}{\noindent\textsc{Proof of Lemma
    \ref{normi}: }}{\hfill$\qedbox$}
\newenvironment{proofcontactpoint}{\noindent\textsc{Proof of Lemma \ref{contactpoint}: }}{\hfill$\qedbox$}
\newenvironment{proofnormpsi}{\noindent\textsc{Proof of Lemma
    \ref{normpsi}: }}{\hfill$\qedbox$}
\newenvironment{prooftheodiffeo}{\noindent\textsc{Proof of Theorem
    \ref{diffeo}: }}{\hfill$\qedbox$}
\newenvironment{proofnirmos}{\noindent\textsc{Proof of Proposition
    \ref{nirmos}: }}{\hfill$\qedbox$}

\newcommand{\scal}[2]{\left\langle #1 , #2\right\rangle}

\def\ndt{\noindent} 
\def\saut{\vspace{0.25cm}} 
\def\R{\mathbb{R}}
\def\S{\mathbb{S}}
\def\insm{\int_{M}} 
\def\vol{dv} 
\def\var{(M^n,g)}
\def\la{\lambda_1(M)}
\def\const{\sqrt{\frac{n}{\la}}}
\def\bonei{(\ei)_{1\leq i\leq n}}
\def\ei{e_i} 
\def\ej{e_j}

\def\dri{\partial_i}
\def\nbo{\nabla^0}
\def\nba{\nabla}

\def\nh{\|H\|_{\infty}}
\def\bh{\|B\|_{\infty}}
\def\nx{\|X\|_{\infty}}
\def\nxti{\|\psi\|_{\infty}}
\def\nfi1{\|\varphi\|_1}
\def\nfii{\|\varphi\|_{\infty}}
\def\nfid{\|\varphi\|_2}
\def\npsi{\|\psi\|_{\infty}}
\def\nfidkmd{\|\xi\|_{2k-2}}
\def\nfiqipi{\|\xi\|_{2q^{i+1}}}
\def\nfiqi{\|\xi\|_{2q^i}^{1-\coefi}}
\def\coefi{\frac{1}{q^i+1}}
\def\coefk{\frac{1}{q^k+1}}
\def\vlph{n\|H\|_{2p}^2}
\def\cepsi{C_{\varepsilon}}

\def\crit{\frac{n}{\la}}
\def\ray{\sqrt{\crit}}
\def\nph2{\|H\|_{2p}^2}

\def\ktld{\tilde{K}}
\def\xt{X^{T}}
\def\nxii{\|\xi\|_{\infty}}
\def\nxid{\|\xi\|_2}
\def\ric{\text{Ric\ }}

\begin{document} 
 
\title {A PINCHING THEOREM FOR THE FIRST EIGENVALUE OF THE LAPLACIAN ON
  HYPERSURFACES OF THE EUCLIDEAN SPACE}

\author{Bruno COLBOIS, Jean-Fran\c{c}ois GROSJEAN} 
\date{ } 
\maketitle

\ndt Supported by European Commission through Human Potential Programme.
\saut

\ndt COLBOIS Bruno, Institut de Math\'ematiques, Universit\'e de Neuch\^atel, Rue \'Emile Argand 11, CH-2007 NEUCH\^ATEL, SUISSE.

\saut

\ndt GROSJEAN Jean-Fran\c cois, Institut \'Elie Cartan (Math\'ematiques), Universit\'e Henri Poinca\-r\'e Nancy I, B.P. 239, F-54506 VANDOEUVRE-LES-NANCY CEDEX, FRANCE. 

\saut 

\ndt E-mail: bruno.colbois@unine.ch, grosjean@iecn.u-nancy.fr
\saut 
 
\begin{abstract} In this paper, we give pinching Theorems for the first nonzero
  eigenvalue $\la$ of the Laplacian on the compact hypersurfaces of
  the Euclidean space. Indeed, we prove that if the volume of $M$ is $1$
  then, for any $\varepsilon>0$, there
  exists a constant $C_{\varepsilon}$ depending on the dimension $n$ of $M$
  and the $L_{\infty}$-norm of the mean curvature $H$,
  so that if the $L_{2p}$-norm $\|H\|_{2p}$ ($p\geq 2$) of $H$ satisfies $\vlph-\cepsi<\la$, then the Hausdorff-distance
  between $M$ and a round sphere of radius $(n/\la)^{1/2}$ is smaller
  than $\varepsilon$. Furthermore, we prove that if
  $C$ is a small enough constant depending on $n$ and the $L_{\infty}$-norm of the
  second fundamental form, then the pinching condition $\vlph-C<\la$ 
  implies that $M$ is diffeomorphic to an $n$-dimensional sphere.

\end{abstract}
\saut

\ndt {\em Key words:} Spectrum, Laplacian, pinching results, 
hypersurfaces.
\saut

\ndt {\em Mathematics Subject Classification (2000):} 53A07, 53C21.

\newpage 
\ \ 
 
\newpage 
 
\section{Introduction and preliminaries} 
\ \ 
Let $\var$ be a compact, connected and oriented $n$-dimensional Riemannian
ma\-nifold without boundary isometrically immersed by $\phi$ into the
$n+1$-dimensional euclidean space $(\R^{n+1}, can)$
(i.e. $\phi^{\star}can=g$). A well known inequality due to Reilly
(\cite{rei}) gives an extrinsic upper bound for the first nonzero eigenvalue $\la$ of the Laplacian of $\var$ in terms of the square of the length of the mean curvature. Indeed, we have

\begin{align}\label{a}\la\leq\frac{n}{V(M)}\insm |H|^2\vol\end{align}

\ndt where $\vol$ and $V(M)$ denote respectively the Riemannian volume
element and the volume of $\var$. Moreover the equality holds if and only
if $\var$ is a geodesic hypersphere of $\R^{n+1}$.

By using H\"older inequality, we obtain some other similar estimates 
for the $L_{2p}$-norm ($p\geq 1$) with $H$ denoted by $\nph2$

\begin{align}\label{b}\la\leq\frac{n}{V(M)^{1/p}}\|H\|_{2p}^2,\end{align}

\ndt and as for the inequality (\ref{a}), the equality case is characterized
by the geodesic hyperspheres of $\R^{n+1}$.

A first natural question is to know if there exists a pinching result as
the one we state now: does a constant
$C$ depending on minimum geometric invariants exist so that if we have the
pinching condition
\saut

\begin{tabular}{cc}$(P_C)$
  &\hspace{3cm} $\displaystyle\frac{n}{V(M)^{1/p}}\|H\|_{2p}^2-C<\la$\end{tabular}
\saut

\ndt then $M$ is close to a sphere in a certain sense?

Such questions are known for the intrinsic lower bound of
Lichnerowicz-Obata (\cite{lic}) of $\la$ in terms of the lower bound of the
Ricci curvature (see \cite{cro}, \cite{ili}, \cite{pet}).
Other pinching results have been proved for Riemannian manifolds
with positive Ricci curvature, with a pinching condition on the $n+1$-st eigenvalue
(\cite{pet}), the diameter (\cite{esc}, \cite{ili}, \cite{wu}), the volume
or the radius (see for instance \cite{col1} and \cite{col2}). 

For instance,
S. Ilias proved in \cite{ili} that there exists $\varepsilon$ depending on
$n$ and an upper bound of the sectional curvature so that if the Ricci
curvature $Ric$ of $M$ satisfies $Ric\geq n-1$ and
$\la\leq\lambda_1(\S^n)+\varepsilon$, then $M$ is homeomorphic to
$\S^n$.  

In this article, we investigate the case of hypersurfaces where, as 
far as we know, very little is known about pinching and stability results
(see however \cite{shixu1}, \cite{shixu2}).

More precisely, in our paper, the hypothesis made in \cite{ili} that $M$ has a positive Ricci curvature is replaced
by the fact that $M$ is isometrically immersed as a hypersurface in 
$\R^{n+1}$, and the bound on the sectional curvature by an $L^{\infty}$-bound on the mean curvature or on the second fundamental form. Note 
that we do not know if such bounds are sharp, or if a bound on the 
$L^q$-norm (for some $q$) of the mean curvature would be enough.

\medskip
We get the following results

\begin{theo}\label{hausdist} Let $\var$ be a compact, connected and
  oriented $n$-dimensional Riemannian ma\-nifold without boundary
  isometrically immersed by $\phi$ in $\R^{n+1}$. Assume that $V(M)=1$ and let $x_0$ be the center of mass of
  $M$. Then for any $p\geq 2$ and for any $\varepsilon>0$, there
  exists a constant $\cepsi$ depending only on $n$, $\varepsilon>0$ and 
  on the $L_{\infty}$-norm of $H$ so that if

\saut

\begin{tabular}{cc}$(P_{\cepsi})$
  &\hspace{3cm} $n\|H\|_{2p}^2-\cepsi<\la$\end{tabular}
\saut

\ndt then the Hausdorff-distance $d_H$ of $M$ to the sphere
$S\left(x_0,\sqrt{\frac{n}{\la}}\right)$ of center $x_0$ and radius
$\sqrt{\frac{n}{\la}}$ satisfies $d_H\left(\phi(M),S\left(x_0,\sqrt{\frac{n}{\la}}\right)\right)<\varepsilon$.

\end{theo}

We recall that the Hausdorff-distance between two compact subsets $A$ and
$B$ of a metric space is given by

$$d_H(A,B)=\inf\{\eta| V_{\eta}(A)\supset B\; \text{and}\;  V_{\eta}(B)\supset
A\}$$

\ndt where for any subset $A$, $V_{\eta}(A)$ is the tubular neighborhood of
$A$ defined by $V_{\eta}(A)= \{x|dist(x, A)<\eta\}$.
\saut

\ndt {\bf Remark} We will see in the proof that $\cepsi(n,\nh)\rightarrow
0$ when $\nh\rightarrow\infty$ or $\varepsilon\rightarrow 0$.
\saut

In fact the previous Theorem is a consequence of the above definition and
the following Theorem

\begin{theo}\label{maintheo} Let $\var$ be a compact, connected and
  oriented $n$-dimensional Riemannian ma\-nifold without boundary
  isometrically immersed by $\phi$ in $\R^{n+1}$. Assume that $V(M)=1$ and
  let $x_0$ be the center of mass of $M$. Then for any $p\geq 2$ and for any $\varepsilon>0$,
  there exists a constant $\cepsi$ depending only on
  $n$, $\varepsilon>0$ and on the $L_{\infty}$-norm of $H$ so that if

\saut

\begin{tabular}{cc}$(P_{\cepsi})$
  &\hspace{3cm} $n\|H\|_{2p}^2-\cepsi<\la$\end{tabular}
\saut

\ndt then

\begin{enumerate} \item $\phi(M)\subset
  B\left(x_0,\const+\varepsilon\right)\backslash
  B\left(x_0,\const-\varepsilon\right)$.

\item $\forall x\in S\left(x_0,\const\right), B(x,\varepsilon)\cap
  \phi(M)\neq\O$.

\end{enumerate}

\end{theo}

In the following Theorem, if the pinching is strong enough, with a 
control on $n$ and the
$L_{\infty}$-norm of the second fundamental form, we obtain 
that $M$ is diffeomorphic to a sphere and even almost isometric with 
a round sphere in a sense we will make precise.

\begin{theo}\label{diffeo} Let $\var$ be a compact, connected and
  oriented $n$-dimensional Riemannian ma\-nifold ($n\geq 2$) without boundary
  isometrically immersed by $\phi$ in $\R^{n+1}$. Assume that
  $V(M)=1$. Then for any $p\geq 2$, there exists a constant $C$ depending
  only on $n$ and  the $L_{\infty}$-norm of the second fundamental form $B$ so that if

\saut

\begin{tabular}{cc}$(P_{C})$
  &\hspace{3cm} $n\|H\|_{2p}^2-C <\la$\end{tabular}
\saut

\ndt Then $M$ is diffeomorphic to $\S^n$.

More precisely, there exists a diffeomorphism $F$ from $M$ into the sphere
$\S^n\left(\ray\right)$ of radius $\ray$ which is a quasi-isometry. 
Namely, for any $\theta$, $0<\theta<1$, there
exists a constant $C$ depending only on $n$, the $L_{\infty}$-norm of $B$
and $\theta$, so that the pinching condition $(P_{C})$ implies

$$\left||dF_x(u)|^2-1\right|\leq \theta$$

\ndt for any $x\in M$ and $u\in T_x M$ so that $|u|=1$.

\end{theo}

Now we will give some preliminaries for the proof of these
Theorems. Throughout the paper, we consider a compact, connected and
oriented $n$-dimensional Riemannian ma\-nifold $\var$ without boundary
isometrically immersed by $\phi$ into $(\R^{n+1}, can)$
(i.e. $\phi^{\star}can=g$). Let $\nu$ be the outward normal vector field. Then the second fundamental form of the immersion will be defined by $B(X,Y)=\scal{\nbo_X \nu}{Y}$,
where $\nbo$ and $\langle \ \ ,\ \ \rangle$ are respectively the Riemannian
connection and the inner product of $\R^{n+1}$. Moreover the mean curvature $H$ will be given by $H=(1/n)trace(B)$.

Now let $\dri$ be an orthonormal frame of $\R^{n+1}$ and let $x_i : \R^{n+1}\rightarrow \R$ be the associated component functions. Putting $X_i= x_i\circ\phi$, a straightforward calculation shows us that

$$B\otimes\nu=-\sum_{i\leq n+1}\nba dX_i\otimes\dri$$

\ndt and

$$nH\nu=\sum_{i\leq n+1}\Delta X_i\dri$$

\ndt where $\nba$ and $\Delta$ denote respectively the Riemannian
connection and the Laplace-Beltrami operator of $\var$. On the other hand,
we have the well known formula

\begin{align}\label{hsiung}\frac{1}{2}\Delta |X|^2=nH\scal{\nu}{X}-n\end{align}

\ndt where $X$ is the position vector given by $X=\sum_{i\leq n+1}X_i\dri$.

We recall that to prove the Reilly inequality, we use the functions $X_i$
as test functions (cf \cite{rei}). Indeed, doing a translation if
necessary, we can assume that $\insm X_i\vol=0$ for all $i\leq n+1$ and we
can apply the variational characterization of $\la$ to $X_i$. If the equality
holds in (\ref{a}) or (\ref{b}), then the functions are nothing but
eigenfunctions of $\la$ and from the Takahashi's Theorem (\cite{tak}) $M$
is immersed isometrically in $\R^{n+1}$ as a geodesic sphere of radius $\sqrt{\crit}$.

Throughout the paper we use some notations. From now on, the inner product
and the norm induced by $g$ and $can$ on a tensor $T$ will be denoted
respectively by $\langle \; ,\;\rangle$ and $|\;|^2$, and the $L_p$-norm will be given by

$$\|T\|_p=\left(\insm|T|^p\vol\right)^{1/p}$$

\ndt and

$$\|T\|_{\infty}=\sup_M |T|$$

We end these preliminaries by a convenient result

\begin{lem}\label{pratic}  Let $\var$ be a compact, connected and
  oriented $n$-dimensional Riemannian ma\-nifold ($n\geq 2$) without boundary
  isometrically immersed by $\phi$ in $\R^{n+1}$. Assume that
  $V(M)=1$. Then there exist constants $c_n$ and $d_n$ depending only on $n$ so
  that for any $p\geq 2$, if $(P_C)$ is true with $C< c_n$ then

\begin{align}\label{titi}\frac{n}{\la}\leq d_n\end{align}

\end{lem}

\saut

\begin{proof} We recall the standard Sobolev inequality (cf \cite{hofspr},
  \cite{hofspr'}, \cite{xu} and p 216 in \cite{burzal}). If $f$ is a smooth function and $f\geq 0$, then

\begin{align}\label{sobol}\left(\insm
    f^{\frac{n}{n-1}}\vol\right)^{1-(1/n)}\leq
  K(n)\insm\left(|df|+|H|f\right)\vol\end{align}

\ndt where $K(n)$ is a constant depending on $n$ and the volume of the unit
ball in $\R^n$. Taking $f=1$ on $M$, and using the fact that $V(M)=1$, we
deduce that 

\begin{align*}\|H\|_{2p}\geq\frac{1}{K(n)}\end{align*}

\ndt and if $(P_C)$ is satisfied and $C\leq\frac{n}{2K(n)^2}=c_n$, then

\begin{align*}\frac{n}{\la}\leq\frac{1}{\vlph-C}\leq 2K(n)^2=d_n\end{align*}

\end{proof}

\saut

Throughout the paper, we will assume that $V(M)=1$ and $\insm X_i\vol=0$
for all $i\leq n+1$. The last assertion implies that the center of mass of
$M$ is the origin of $\R^{n+1}$.


\section{An $L^2$-approach of the problem}
\ \

A first step in the proof of the Theorem \ref{maintheo} is to
prove that if the pinching condition $(P_C)$ is satisfied, then $M$ is close
to a sphere in an $L^2$-sense.

In the following Lemma, we prove that the $L^2$-norm of the position vector
is close to $\sqrt\crit$.

\begin{lem}\label{vectpos} If we have the pinching condition $(P_C)$ with $C<c_n$, then

$$\frac{n\la}{(C+\la)^2}\leq\|X\|_2^2\leq \frac{n}{\la}\leq d_n$$

\end{lem}

\begin{proof} Since $\insm X_i\vol=0$, we can apply the variational
characterization of the eigenvalues to obtain

\begin{align*}\la\insm\sum_{i\leq
    n+1}|X_i|^2\vol\leq\insm\sum_{i\leq n+1}|dX_i|^2\vol=n\end{align*}

\ndt which gives the inequality of the right-hand side

Let us prove now the inequality of the left-hand side.

\begin{align*}\la\insm|X|^2\vol&\leq\frac{\left(\insm\displaystyle\sum_{i\leq n+1}|dX_i|^2\vol\right)^4}{\left(\insm\displaystyle\sum_{i\leq n+1}|dX_i|^2\vol\right)^3}=\frac{\left(\insm\displaystyle\sum_{i\leq n+1}(\Delta X_i) X_i\vol\right)^4}{n^3}\displaybreak[2]\\
&\leq\frac{\left(\insm\displaystyle \sum_{i\leq n+1}(\Delta
    X_i)^2\vol\right)^2\left(\insm|X|^2\vol\right)^2}{n^3}\displaybreak[2]\\
    &= n\left(\insm
    H^2\vol\right)^2\left(\insm|X|^2\vol\right)^2\end{align*}

\ndt then using again the H\"older inequality, we get

\begin{align*}\la&\leq\frac{1}{n}\left(\vlph\right)^2\insm|X|^2\vol\leq\frac{(C+\la)^2}{n}\insm|X|^2\vol\end{align*}

\ndt This completes the proof. \end{proof}
\saut

From now on, we will denote by $\xt$ the orthogonal tangential 
  projection on $M$. In fact, at $x\in M$, $\xt$ is nothing but the 
  vector of $T_{x}M$ defined by $\xt=\displaystyle\sum_{1\leq i\leq 
  n}\scal{X}{\ei}\ei$ where $\bonei$ is an orthonormal basis of $T_{x}M$.
In the following Lemma, we will show that the condition $(P_C)$ implies
that the $L^2$-norm of $\xt$ of $X$ on $M$ is close to $0$. 

\begin{lem}\label{l2xt} If we have the pinching condition $(P_C)$, then

$$\|\xt\|_2^2\leq A(n)C$$

\end{lem}

\begin{proof} From the lemma \ref{vectpos} and the relation (\ref{hsiung}), we have

\begin{align*}\la\insm|X|^2\vol&\leq n=n\left(\insm H\scal{X}{\nu}\vol\right)^2\displaybreak[2]\\
&\leq\left(\insm |H||\scal{X}{\nu}|\vol\right)^2\leq\vlph\left(\insm
  |\scal{X}{\nu}|^{\frac{2p}{2p-1}}\vol\right)^{\frac{2p-1}{p}}\displaybreak[2]\\
&\leq\vlph\left(\insm
  |\scal{X}{\nu}|^2\vol\right)=\vlph\insm|X|^2\vol\end{align*}

\ndt Then we deduce that

\begin{align*}\vlph\|\xt\|_2^2&=\vlph\left(\insm\left(|X|^2-|\scal{X}{\nu}|^2\right)\vol\right)\displaybreak[2]\\
&\leq (\vlph-\la)\|X\|_2^2\leq d_n C\end{align*}

\ndt where in the last inequality we have used the pinching condition and
the Lemma \ref{vectpos}.\end{proof}

\saut

Now, we will show that the condition $(P_C)$ implies that
the component functions are almost eigenfunctions in an $L^2$-sense. For
this, let us consider the vector field $Y$ on $M$ defined by

$$Y=\sum_{i\leq n+1}\left(\Delta X_i-\la X_i\right)\dri=nH\nu-\la X$$

\begin{lem}\label{Y} If $(P_C)$ is satisfied, then

$$\|Y\|_2^2\leq nC$$ 

\end{lem}

\saut

\begin{proof} We have

\begin{align*}\insm |Y|^2\vol=\insm\left(n^2H^2-2n\la H\scal{\nu}{X}+\la^2|X|^2\right)\vol\end{align*}

\ndt Now by integrating the relation (\ref{hsiung}) we deduce that

$$\insm H\scal{\nu}{X}\vol=1$$

\ndt Furthermore, since $\insm X_i\vol=0$, we can apply the variational
characterization of the eigenvalues to obtain

\begin{align*}\la\insm|X|^2\vol=\la\insm\sum_{i\leq n+1}|X_i|^2\vol\leq\insm\sum_{i\leq n+1}|dX_i|^2\vol=n\end{align*}

\ndt Then

\begin{align*}\insm |Y|^2\vol&\leq n^2\insm|H|^2\vol-n\la \leq n \left(\vlph-\la\right)\leq nC\end{align*}

\ndt where in this last inequality we have used the H\"older inequality. \end{proof}
\saut

To prove Assertion 1 of Theorem \ref{maintheo}, we will show that
$\left\| |X|-\left(\crit\right)^{1/2}\right\|_{\infty}\leq\varepsilon$. For this
we need to have an $L^2$-upper bound on the function
$\varphi=|X|\left(|X|-\left(\crit\right)^{1/2}\right)^2$.

Before giving such estimate, we will introduce the vector field $Z$ on $M$
defined by

$$Z=\left(\crit\right)^{1/2}|X|^{1/2}H\nu-\frac{X}{|X|^{1/2}}$$

\ndt We have

\begin{lem}\label{Z} If $(P_C)$ is satisfied with $C<c_n$, then

$$\|Z\|_2^2\leq B(n)C$$

\end{lem}

\begin{proof} We have

\begin{align*}\|Z\|_2^2&=\left\|\left(\crit\right)^{1/2}|X|^{1/2}H\nu-\frac{X}{|X|^{1/2}}\right\|_2^2\displaybreak[2]\\
&=\insm\left(\crit |X|H^2-2\left(\crit\right)^{1/2}H\scal{\nu}{X}+|X|\right)\vol\\
&\leq\crit\left(\insm|X|^2\vol\right)^{1/2}\left(\insm H^4\vol\right)^{1/2}-2\left(\crit\right)^{1/2}+\left(\insm|X|^2\vol\right)^{1/2}\end{align*}

\ndt Note that we have used the relation (\ref{hsiung}). Finally for $p\geq 2$,
we get

\begin{align*}\|Z\|_2^2&\leq\left(\insm|X|^2\vol\right)^{1/2}\left(\crit\nph2+1\right)-2\left(\crit\right)^{1/2}\displaybreak[2]\\
&\leq\left(\crit\right)^{1/2} \left(\frac{C}{\la}+2\right)-2\left(\crit\right)^{1/2}
\\
&=\left(\crit\right)^{1/2}\frac{C}{\la}\leq\frac{d_n^{3/2}}{n}C\end{align*}

\ndt This concludes the proof of the Lemma. \end{proof}

\saut

\ndt Now we give an $L^2$-upper bound of $\varphi$

\saut

\begin{lem}\label{norml2} Let $p\geq 2$ and $C\leq c_n$. If we have the pinching condition $(P_C)$, then

$$\|\varphi\|_2\leq D(n)\nfii^{3/4}C^{1/4}$$

\end{lem}

\saut 

\begin{proof} We have

$$\|\varphi\|_2=\left(\insm\varphi^{3/2}\varphi^{1/2}\vol\right)^{1/2}\leq\nfii^{3/4}\|\varphi^{1/2}\|_1^{1/2}$$

\ndt and noting that

$$|X|\left(|X|-\left(\crit\right)^{1/2}\right)^2=\left|
  |X|^{1/2}X-\left(\crit\right)^{1/2}\frac{X}{|X|^{1/2}}\right|^2$$

\ndt we get

\begin{align}\label{ineq}\insm\varphi^{1/2}\vol&=\left\||X|^{1/2}X-\left(\crit\right)^{1/2}\frac{X}{|X|^{1/2}}\right\|_1\notag\\
&=\left\|-\frac{|X|^{1/2}}{\la}Y+\crit|X|^{1/2}H\nu-\left(\crit\right)^{1/2}\frac{X}{|X|^{1/2}}\right\|_1\notag\displaybreak[2]\\
&\leq\left\|\frac{|X|^{1/2}}{\la}Y\right\|_1+\left(\crit\right)^{1/2}\left\|Z\right\|_1\end{align}

\ndt From Lemmas \ref{Y} and \ref{pratic} we get

\begin{align*}\left\|\frac{|X|^{1/2}}{\la}Y\right\|_1&\leq\frac{1}{\la}\left(\insm
  |X|\vol\right)^{1/2}\|Y\|_2\\
&\leq \frac{1}{\la}\left(\insm
  |X|^2\vol\right)^{1/4}\|Y\|_2\leq\frac{d_n^{3/4}}{n^{1/2}}C^{1/2}\end{align*}

\ndt Moreover, using Lemmas \ref{Z} and \ref{pratic} again it is easy to see that the last term of (\ref{ineq}) is bounded
by $d_n^{1/2}B(n)^{1/2}C^{1/2}$. Then $\|\varphi^{1/2}\|_1^{1/2}\leq D(n)C^{1/4}$.
 
\end{proof}


\section{Proof of Theorem \ref{maintheo}}
\ \

The proof of Theorem \ref{maintheo} is immediate from the two following
technical Lemmas which we state below.

\begin{lem}\label{normi} For $p\geq 2$ and for any $\eta>0$, there exists
  $K_{\eta}(n, \nh)\leq c_n$ so that if $(P_{K_{\eta}})$ is true, then
  $\|\varphi\|_{\infty}\leq\eta$. Moreover, $K_{\eta}\rightarrow 0$ when
  $\nh\rightarrow\infty$ or $\eta\rightarrow 0$.

\end{lem}
\saut

\ndt and

\begin{lem}\label{contactpoint} Let $x_0$ be a point of the sphere $S(O,R)$
  of $\R^{n+1}$ with the center at the origin and of radius $R$. Assume that $x_0=Re$ where $e\in\S^n$. Now let $\var$ be a compact oriented
  $n$-dimensional Riemannian manifold without boundary isometrically
  immersed by $\phi$ in $\R^{n+1}$ so that
  $\phi(M)\subset\left(B(O,R+\eta)\backslash B(O,R-\eta)\right)\backslash
  B(x_0,\rho)$ with $\rho=4(2n-1)\eta$ and suppose that there exists
  a point $p\in M$ so that $\scal{X}{e}>0$. Then there exists $y_0\in M$
  so that the mean curvature $H(y_0)$ at $y_0$ satisfies $|H(y_0)|\geq \frac{1}{4n\eta}$.

\end{lem}

Now, let us see how to use these Lemmas to prove Theorem \ref{maintheo}.
\saut

\begin{proofmaintheo} Let $\varepsilon>0$ and let us consider the function
  $f(t)=t\left(t-\left(\crit\right)^{1/2}\right)^2$. Let us put

\begin{align*}\eta(\varepsilon)&=\min\left(\left(\frac{1}{\nh}-\varepsilon\right)\varepsilon^2,\left(\frac{1}{\nh}+\varepsilon\right)\varepsilon^2,\frac{1}{27\nh^3}\right)\\
&\leq\min\left(f\left(\left(\crit\right)^{1/2}-\varepsilon\right),f\left(\left(\crit\right)^{1/2}+\varepsilon\right),\frac{1}{27\nh^3}\right)\end{align*}

\ndt  Then, as $\eta(\varepsilon)>0$ and from Lemma \ref{normi}, it follows that if the pinching condition
$(P_{K_{\eta(\varepsilon)}})$ is satisfied with $K_{\eta(\varepsilon)}\leq c_n$, then for any $x\in M$, we have

\begin{align}\label{major}f(|X|)\leq\eta(\varepsilon)\end{align}

\ndt Now to prove Theorem \ref{maintheo}, it is sufficient to assume 
$\varepsilon<\frac{2}{3\nh}$. Let us show that either

\begin{align}\label{etude}\left(\crit\right)^{1/2}-\varepsilon\leq |X|\leq
\left(\crit\right)^{1/2}+\varepsilon\;\;\;\;\text{or}\;\;\;\;
|X|<\frac{1}{3}\left(\crit\right)^{1/2}\end{align}

\ndt By studying the function $f$, it is easy to see that $f$ has a
unique local maximum in $\frac{1}{3}\left(\crit\right)^{1/2}$ and from the
definition of $\eta(\varepsilon)$ we have $\eta(\varepsilon)<\frac{4}{27}\frac{1}{\nh^3}\leq\frac{4}{27}\left(\crit\right)^{3/2}=f\left(\frac{1}{3}\left(\crit\right)^{1/2}\right)$.

Now since $\varepsilon<\frac{2}{3\nh}$, we have
$\varepsilon<\frac{2}{3}\left(\crit\right)^{1/2}$, and 
$\frac{1}{3}\left(\crit\right)^{1/2}<\left(\crit\right)^{1/2}-\varepsilon$. 
This and (\ref{major}) yield (\ref{etude}).

 Now, from Lemma
\ref{vectpos} we deduce that there exists a point $y_0\in M$ so that
$|X(y_0)|\geq\frac{n^{1/2}\la^{1/2}}{(K_{\eta(\varepsilon)}+\la)}$ and since
$K_{\eta(\varepsilon)}\leq c_n=\frac{n}{d_n}\leq \la\leq 2\la$ (see the
proof of the Lemma \ref{pratic}), we obtain $|X(y_0)|\geq\frac{1}{3}\left(\crit\right)^{1/2}$.

\ndt By the connectedness of $M$, it follows that $\left(\crit\right)^{1/2}-\varepsilon\leq|X|\leq
\left(\crit\right)^{1/2}+\varepsilon$ for any point of $M$ and
Assertion 1 of Theorem \ref{maintheo} is shown for the condition $(P_{K_{\eta(\varepsilon)}})$.

In order to prove the second assertion, let us consider the pinching condition
$(P_{C_{\varepsilon}})$ with 
$C_{\varepsilon}=K_{\eta\left(\frac{\varepsilon}{4(2n-1)}\right)}$. Then Assertion 1
is still valid. Let $x=\left(\crit\right)^{1/2}e\in
S\left(O,\const\right)$, with $e\in\S^n$ and suppose that
$B(x,\varepsilon)\cap M=\O$. Since $\insm X_i\vol=0$ for any
$i\leq n+1$, there exists a point $p\in M$ so that $\scal{X}{e}>0$ and we
can apply Lemma \ref{contactpoint}. Therefore there is a point $y_0\in
M$ so that $H(y_0)\geq\frac{2n-1}{n\varepsilon}>\nh$ since we have 
assumed $\varepsilon<\frac{2}{3\nh}\leq\frac{2n-1}{2n\nh}$. Then we 
obtain a contradiction which implies $B(x,\varepsilon)\cap M\neq\O$ and
Assertion 2 is satisfied. Furthermore,
$C_{\varepsilon}\rightarrow 0$ when $\nh\rightarrow\infty$ or
$\varepsilon\rightarrow 0$. \end{proofmaintheo}


\section{Proof of Theorem \ref{diffeo}}
\ \

From Theorem \ref{maintheo}, we know that for any $\varepsilon>0$, there
exists $\cepsi$ depending only on $n$ and $\nh$ so that if $(P_{\cepsi})$ is true then

$$\left| |X|_x-\ray\right|\leq\varepsilon$$

\ndt for any $x\in M$. Now, since $\sqrt{n}\nh\leq\bh$, it is easy to see from the previous
proofs that we can assume that $\cepsi$ is depending only on $n$ and
$\bh$.

The proof of Theorem \ref{diffeo} is a consequence of the following
Lemma on the $L_{\infty}$-norm of $\psi=|\xt|$

\begin{lem}\label{normpsi} For $p\geq 2$ and for any $\eta>0$, there exists
  $K_{\eta}(n,\bh)$ so that if $(P_{K_{\eta}})$ is true, then
$\nxti\leq\eta$.  Moreover, $K_{\eta}\rightarrow 0$ when
  $\bh\rightarrow\infty$ or $\eta\rightarrow 0$.

\end{lem}

\ndt This Lemma will be proved in the Section 5.

\saut

\begin{prooftheodiffeo} Let
  $\varepsilon<\frac{1}{2}\sqrt{\frac{n}{\bh}}\leq\ray$. From the choice of
  $\varepsilon$, we deduce that the condition $(P_{\cepsi})$ implies that
  $|X_x|$ is nonzero for any $x\in M$ (see the proof of Theorem \ref{maintheo}) and we can consider the differential application

\begin{center}\begin{tabular}{ccccc }$F$ & $:$ & $M$ & $\longrightarrow$ &
    $S\left(O,\ray\right)$\\
\ \ &\ \ & $x$ & $\longmapsto$ &
$\sqrt{\crit}\frac{X_x}{|X_x|}$\end{tabular}\end{center}

\ndt We will prove that $F$ is a quasi isometry. Indeed, for any
$0<\theta<1$, we can choose a
constant $\varepsilon(n,\bh, \theta)$ so that for any $x\in M$ and any unit
vector $u\in T_x M$, the pinching condition $(P_{C_{\varepsilon(n,\bh,\theta)}})$ implies

$$\left||dF_x(u)|^2-1\right|\leq\theta$$

\ndt For this, let us compute $dF_x(u)$. We have

\begin{align*}dF_x(u)&=\ray\nbo_u\left(\frac{X}{|X|}\right)\Big{|}_x=\ray u\left(\frac{1}{|X|}\right)X+\ray\frac{1}{|X|}\nbo_u X\displaybreak[2]\\
&=-\frac{1}{2}\ray\frac{1}{|X|^3}u(|X|^2)X+\ray\frac{1}{|X|}u\\
&=-\ray\frac{1}{|X|^3}\scal{u}{X}X+\ray\frac{1}{|X|}u\displaybreak[2]\\
&=\ray\frac{1}{|X|}\left(-\frac{\scal{u}{X}}{|X|^2}X+u\right)\end{align*}

\ndt By a straightforward computation, we obtain

\begin{align}\label{quasisom}\left||dF_x(u)|^2-1\right|&=\left|\crit\frac{1}{|X|^2}\left(1-\frac{\scal{u}{X}^2}{|X|^2}\right)-1\right|\notag\\
&\leq\left|\crit\frac{1}{|X|^2}-1\right|+\crit\frac{1}{|X|^4}\scal{u}{X}^2\end{align}

\ndt Now

\begin{align*}\left|\crit\frac{1}{|X|^2}-1\right|=\frac{1}{|X|^2}\left|\crit-|X|^2\right|&\leq\varepsilon\frac{\left|\ray+|X|\right|}{|X|^2}\leq\varepsilon\frac{2\ray+\varepsilon}{\left(\ray-\varepsilon\right)^2}\end{align*}

\ndt Let us recall that $\frac{n}{d_n}\leq \la\leq \bh^2$ (see (\ref{titi}) for the first
inequality). Since we assume $\varepsilon<\frac{1}{2}\sqrt{\frac{n}{\bh}}$, the right-hand side is bounded above by a constant depending only 
on $n$ and $\bh$ and we have

\begin{align}\label{pinc1}\left|\crit\frac{1}{|X|^2}-1\right|\leq\varepsilon\gamma(n,\bh)\end{align}

\ndt On the other hand, since $\cepsi(n,\bh)\rightarrow 0$ when
$\varepsilon\rightarrow 0$, there exists $\varepsilon(n,\bh,\eta)$ so that
$C_{\varepsilon_{(n,\bh,\eta)}}\leq K_{\eta}(n,\bh)$ (where $K_{\eta}$ is
the constant of the Lemma) and then by Lemma
\ref{normpsi}, $\nxti^2\leq\eta^2$. 
Thus, there exists a constant $\delta$ depending only on $n$ and 
$\bh$ so that

\begin{align}\label{pinc2}\crit\frac{1}{|X|^4}\scal{u}{X}^2\leq\crit\frac{1}{|X|^4}\nxti^2\leq\eta^2\delta(n,\bh)\end{align}

\ndt then from (\ref{quasisom}), (\ref{pinc1}) and (\ref{pinc2}) we deduce that the
condition $(P_{C_{\varepsilon(n,\bh,\eta)}})$ implies

$$\left||dF_x(u)|^2-1\right|\leq\varepsilon\gamma(n,\bh)+\eta^2\delta(n,\bh)$$

\ndt Now let us choose
$\eta=\left(\frac{\theta}{2\delta}\right)^{1/2}$. Then we can assume that
$\varepsilon(n,\bh,\eta)$ is small enough in order to have
$\varepsilon(n,\bh,\eta)\gamma(n\bh)\leq\frac{\theta}{2}$. In this case we have

$$\left||dF_x(u)|^2-1\right|\leq \theta$$

\ndt Now let us fix $\theta$, $0<\theta<1$. It follows that $F$ is a local
diffeomorphism from $M$ to $S\left(O,\ray\right)$. Since
$S\left(O,\ray\right)$ is simply connected for $n\geq 2$, $F$ is a diffeomorphism. \end{prooftheodiffeo}


\section{Proof of the technical Lemmas }
\ \

The proofs of Lemmas \ref{normi} and \ref{normpsi} are providing from a
result stated in the following Proposition using a Nirenberg-Moser type of proof.

\begin{pro}\label{nirmos} Let $\var$ be a compact, connected and oriented $n$-dimensional Riemannian
ma\-nifold without boundary isometrically immersed into the
$n+1$-dimensional euclidean space $(\R^{n+1}, can)$. Let $\xi$ be a
nonnegative continuous function so that $\xi^k$ is smooth for $k\geq
2$. Let $0\leq r<s\leq 2$ so that

$$\frac{1}{2}\Delta\xi^2
\xi^{2k-2}\leq\delta\omega+(A_1+kA_2)\xi^{2k-r}+(B_1+kB_2)\xi^{2k-s}$$

\ndt where $\delta\omega$ is the codifferential of a $1$-form and $A_1,A_2,B_1,B_2$ are nonnegative constants. Then for any
$\eta>0$, there exists a constant $L(n,A_1,A_2,B_1,B_2,\nh,\eta)$ depending only on $n$,
$A_1$, $A_2$, $B_1$, $B_2$, $\nh$ and $\eta$ so that if $\nxii>\eta$ then

$$\nxii\leq L(n,A_1, A_2, B_1, B_2, \nh,\eta)\nxid$$

\ndt Moreover, $L$ is bounded when $\eta\rightarrow\infty$, and if $B_1>0$, $L\rightarrow\infty$ when $\nh\rightarrow\infty$ or
$\eta\rightarrow 0$.

\end{pro}

\saut

This Proposition will be proved at the end of the paper.
\saut

Before giving the proofs of Lemmas \ref{normi} and \ref{normpsi}, we
will show that under the pinching condition $(P_C)$ with $C$ small enough,
the $L_{\infty}$-norm of $X$ is bounded by a constant depending only on $n$ and $\nh$.

\begin{lem}\label{ninfx} If we have the pinching condition $(P_C)$ with
  $C<c_n$, then there exists $E(n,\nh)$ depending only on $n$ and $\nh$ so
  that $\nx\leq E(n,\nh)$.

\end{lem}

\begin{proof} From the relation (\ref{hsiung}), we have

$$\frac{1}{2}\Delta |X|^2 |X|^{2k-2}\leq n\nh |X|^{2k-1}$$

\ndt Then applying Proposition \ref{nirmos} to the function $\xi=|X|$
with $r=0$ and $s=1$, we obtain that if $\nx > E$, then there
exists a constant $L(n,\nh,E)$ depending only on $n$, $\nh$ and $E$
so that

$$\nx\leq L(n,\nh,E)\|X\|_2$$

\ndt and under the pinching condition  $(P_C)$ with $C<c_n$ we have from
Lemma \ref{vectpos}

$$\nx\leq L(n,\nh,E)d_n^{1/2}$$

\ndt Now since $L$ is bounded when
$E\rightarrow\infty$, we can choose $E=E(n,\nh)$ great enough so that

$$L(n,\nh,E)d_n^{1/2}<E$$

\ndt In this case, we have $\nx\leq E(n,\nh)$. 

\end{proof}

\saut

\begin{proofnormi} First we compute the Laplacian of the square of
  $\varphi^2$. We have

\begin{align*}
  \Delta\varphi^2&=\Delta\left(|X|^4-2\left(\crit\right)^{1/2}|X|^3+\crit
    |X|^2\right)\\
    &=-2|X|^2|d|X|^2|^2+2|X|^2\Delta|X|^2\\
&-2\left(\crit\right)^{1/2}\left(-\frac{3}{4}|X|^{-1}|d|X|^2|^2+\frac{3}{2}|X|\Delta|X|^2\right)+\crit\Delta|X|^2\end{align*}

\ndt Now by a direct computation one gets $|d|X|^2|^2\leq 4|X|^2$. Moreover
by the relation (\ref{hsiung}) we have $|\Delta|X|^2|\leq 2n\nh |X|+n$. Then
applying Lemmas \ref{pratic} and \ref{ninfx} we get

$$\Delta\varphi^2\leq \alpha(n,\nh)$$

\ndt and

$$\frac{1}{2}\Delta\varphi^2\varphi^{2k-2}\leq
\alpha(n,\nh)\varphi^{2k-2}$$

\ndt Now, we apply Proposition \ref{nirmos} with $r=0$ and $s=2$. Then
if $\nfii>\eta$, there exists a constant $L(n,\nh)$ depending only on
$n$ and $\nh$ so that

$$\nfii\leq L\nfid$$

\ndt From Lemma \ref{norml2}, if $C\leq c_n$ and
$(P_C)$ is true, we have $\nfid\leq D(n)\nfii^{3/4}C^{1/4}$. Therefore

$$\nfii\leq (L D)^4C$$

\ndt Consequently, if we choose $C=K_{\eta}=\inf\left(\frac{\eta}{(LD)^4},c_n\right)$, then we obtain that $\nfii\leq\eta$.

\end{proofnormi}

\saut

\begin{proofnormpsi} First we will prove that for any $C<c_n$, if $(P_C)$ is
  true, then

\begin{align}\label{deltapsi}\frac{1}{2}(\Delta\psi^2)\psi^{2k-2}\leq\delta\omega+\left( \alpha_1 (n,\bh)+k\alpha_2 (n,\bh) \right) \psi^{2k-2}\end{align}

\ndt where $\delta\omega$ is the codifferential of a $1$-form $\omega$.

First observe that the gradient $\nabla^M |X|^2$ of $|X|^2$ satisfies $\nabla^M |X|^2=2\xt$. Then by the Bochner formula we get

\begin{align*}\frac{1}{2}\Delta|\xt|^2&=\frac{1}{4}\scal{\Delta d|X|^2}{d|X|^2}-\frac{1}{4}|\nabla d|X|^2|^2-\frac{1}{4}\ric (\nabla^M |X|^2,\nabla^M |X|^2)\notag\\
&\leq \frac{1}{4}\scal{d\Delta |X|^2}{d|X|^2}-\frac{1}{4}\ric (\nabla^M |X|^2,\nabla^M |X|^2)\end{align*}

\ndt and by the Gauss formula we obtain

\begin{align*}\frac{1}{2}\Delta|\xt|^2&\leq \frac{1}{4}\scal{d\Delta |X|^2}{d|X|^2}-\frac{1}{4}nH\scal{B\nabla^M |X|^2}{\nabla^M |X|^2}+\frac{1}{4}|B\nabla^M |X|^2|^2\\
&=\frac{1}{4}\scal{d\Delta |X|^2}{d|X|^2}-nH\scal{B\xt}{\xt}+|B\xt|^2
\end{align*}

\ndt By Lemma \ref{ninfx} we know that $\nx\leq E(n,\bh)$ (the
dependance in $\nh$ can be replaced by $\bh$). Then it follows that

\begin{align}\label{estipsi}\frac{1}{2}(\Delta\psi^2)\psi^{2k-2}\leq \frac{1}{4}\scal{d\Delta |X|^2}{d|X|^2}\psi^{2k-2}+\alpha'(n,\bh)\psi^{2k-2}\end{align}

\ndt Now, let us compute the term $\scal{d\Delta |X|^2}{d|X|^2}\psi^{2k-2}$. We
have

\begin{align*}\scal{d\Delta |X|^2}{d|X|^2}\psi^{2k-2}&=\delta\omega+(\Delta |X|^2)^2\psi^{2k-2}-(2k-2)\Delta |X|^2\scal{d|X|^2}{d\psi}\psi^{2k-3}\\
&=\delta\omega+(\Delta |X|^2)^2\psi^{2k-2}-2(2k-2)\Delta |X|^2\scal{\xt}{\nabla^M \psi}\psi^{2k-3}\end{align*}

\ndt where $\omega=-\Delta |X|^2 \psi^{2k-2} d|X|^2$. Now,

$$\ei(\psi)=\frac{\ei|\xt|^2}{2|\xt|}=\frac{\ei|X|^2-\ei\scal{X}{\nu}^2}{2|\xt|}=\frac{\scal{\ei}{X}-B_{ij}\scal{X}{\ej}\scal{X}{\nu}}{|\xt|}$$

\ndt Then

\begin{align*}\scal{d\Delta |X|^2}{d|X|^2}\psi^{2k-2}&=\delta\omega+(\Delta |X|^2)^2\psi^{2k-2}-2(2k-2)\Delta |X|^2 |\xt|\psi^{2k-3}\\
&\hspace{4cm}+2(2k-2)\Delta |X|^2\frac{\scal{B\xt}{\xt}}{|\xt|}\scal{X}{\nu}\psi^{2k-3}\\
&\leq \delta\omega+(\Delta |X|^2)^2\psi^{2k-2}+2(2k-2)|\Delta |X|^2|\psi^{2k-2}\\
&\hspace{4cm}+2(2k-2)|\Delta |X|^2||B||X|\psi^{2k-2}\end{align*}

\ndt Now by relation (\ref{hsiung}) and Lemma \ref{ninfx} we have

$$\scal{d\Delta |X|^2}{d|X|^2}\psi^{2k-2}\leq \delta\omega+\left( \alpha_1''(n,\bh)+k\alpha_2''(n,\bh)\right)\psi^{2k-2}$$

\ndt Inserting this in (\ref{estipsi}), we obtain the desired
inequality (\ref{deltapsi}).

Now applying again Proposition \ref{nirmos}, we get that there exists
$L(n,\bh,\eta)$ so that if $\npsi>\eta$ then

$$\npsi\leq L \|\psi\|_2$$

\ndt From the Lemma \ref{l2xt} we deduce that if the pinching condition
$(P_C)$ holds then $\|\psi\|_2\leq A(n)^{1/2}C^{1/2}$. Then taking
$C=K_{\eta}=\inf\left(\frac{\eta}{L A^{1/2}},c_n\right)$, then
$\npsi\leq\eta$.

\end{proofnormpsi}

\saut

\begin{proofcontactpoint} The idea of the proof consists in 
foliating the region $B(O,R+\eta) \backslash B(O,R-\eta)$ with 
hypersurfaces of large mean curvature and to show that one of these 
hypersurfaces is tangent to $\phi(M)$. This will imply that $\phi(M)$ 
has a large mean curvature at the contact point. 

Consider $\S^{n-1}\subset\R^n$ and
  $\R^{n+1}=\R^n\times\R e$. Let $a$, $L>l>0$ and

\begin{center}\begin{tabular}{ccccc}$\Phi_{L,l,a}$ & $:$ & 
$\S^{n-1}\times\S^1$ & $\longrightarrow$ & $\R^{n+1}$\\
\ \ &\ \ & $(\xi,\theta)$ & $\longmapsto$ &
$L\xi-l\cos\theta\xi+l\sin\theta e+ae$\end{tabular}\end{center}

\ndt Then $\Phi_{L,l,a}$ is a family of embeddings from 
$\S^{n-1}\times\S^1$ in $\R^{n+1}$. If we orient the family 
of hypersurfaces $\Phi_{L,l,a}(\S^{n-1}\times\S^1)$ by the unit outward
normal vector field, a straightforward computation shows
that the mean curvature $H(\theta)$ depends only on $\theta$ and we
have

\begin{align}\label{meancurv}H(\theta)=\frac{1}{n}\left(\frac{1}{l}-\frac{(n-1)\cos\theta}{L-l\cos\theta}\right)\geq\frac{1}{n}\left(\frac{1}{l}-\frac{n-1}{L-l}\right)\end{align}

\ndt Now, let us consider the hypotheses of the Lemma and for
$t_0=2\arcsin\left(\frac{\rho}{2R}\right)\leq t\leq\frac{\pi}{2}$, put
$L=R\sin t$, $l=2\eta$ and $a=R\cos t$. Then $L>l$ and we can consider for $t_0\leq
t\leq\frac{\pi}{2}$ the family ${\cal M}_{R,\eta, t}$ of hypersurfaces
defined by ${\cal M}_{R,\eta, t}=\Phi_{R\sin t, 2\eta, R\cos 
t}(\S^{n-1}\times\S^1)$.

From the relation (\ref{meancurv}), the mean curvature $H_{R,\eta, t}$ of
${\cal M}_{R,\eta, t}$ satisfies

\begin{align*}H_{R,\eta,
    t}&\geq\frac{1}{n}\left(\frac{1}{2\eta}-\frac{n-1}{R\sin t
      -2\eta}\right)\geq\frac{1}{n}\left(\frac{1}{2\eta}-\frac{n-1}{R\sin t_0
      -2\eta}\right)\displaybreak[2]\\
&\geq\frac{1}{n}\left(\frac{1}{2\eta}-\frac{n-1}{R\sin(t_0/2)
      -2\eta}\right)=\frac{1}{n}\left(\frac{1}{2\eta}-\frac{n-1}{\frac{\rho}{2}
      -2\eta}\right)=\frac{1}{4n\eta}\displaybreak[2]\end{align*}

\ndt where we have used in this last equality the fact that
$\rho=4(2n-1)\eta$.

Since there exists a point $p\in M$ so that $\scal{X(p)}{e}>0$, we can find
$t\in [t_0, \pi/2]$ and a point $y_0\in M$ which is a contact point with
${\cal M}_{R,\eta, t}$. Therefore $|H(y_0)|\geq \frac{1}{4n\eta}$.

\includegraphics[height=10cm]{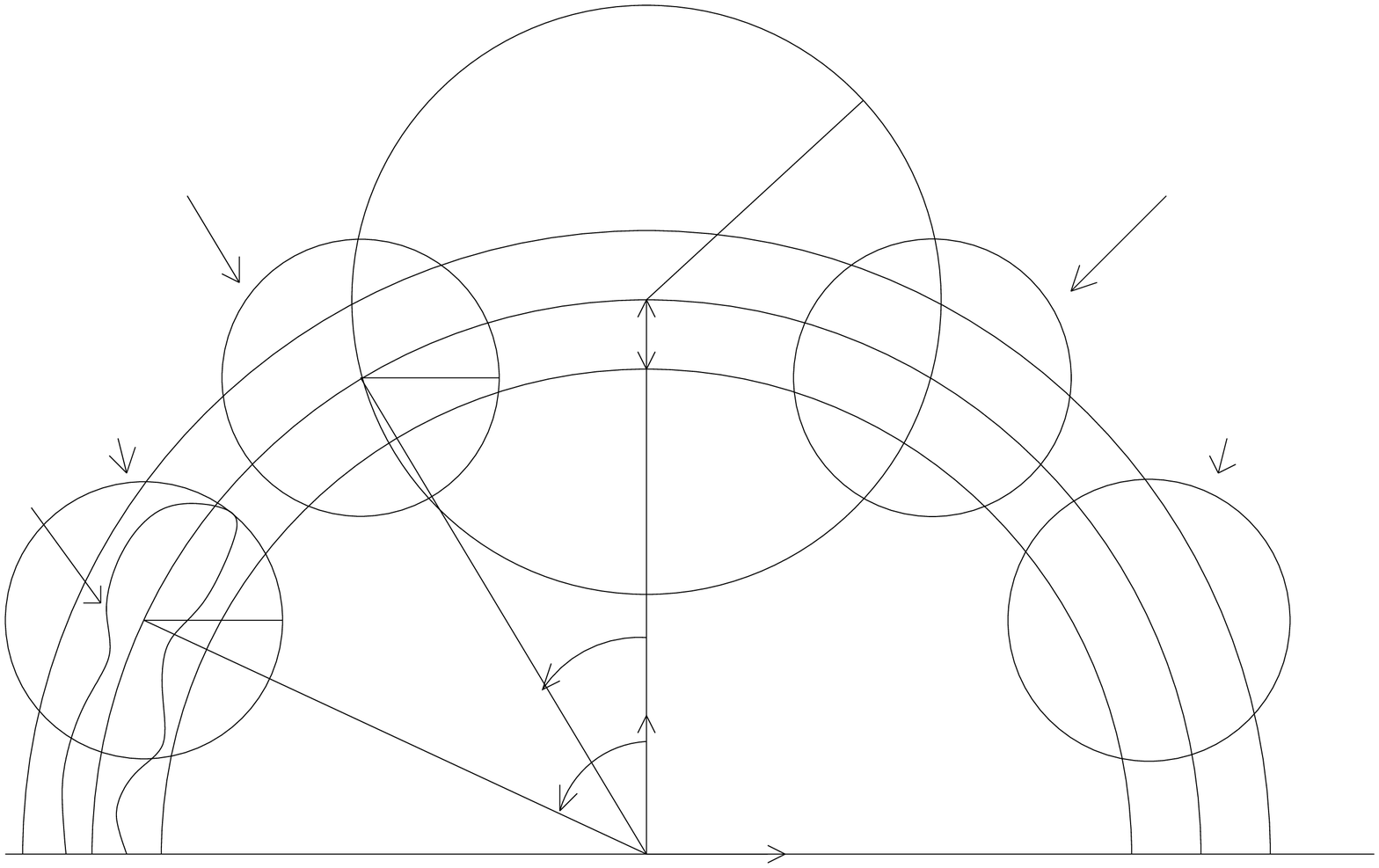}
\vspace{-8,9cm}

\hspace{0,5cm} {\scriptsize ${\cal M}_{R,\eta, t_0}\cap F$ } \hspace{5,8cm}
{\scriptsize $\rho$} \hspace{3cm} {\scriptsize ${\cal M}_{R,\eta, t_0}\cap F$ }

\vspace{0,7cm}

\hspace{7cm}{\scriptsize $x_0$}

\vspace{0,18cm}

\hspace{7,5cm}{\scriptsize $\eta$}

\vspace{0,2cm}

\hspace{0,5cm} {\scriptsize ${\cal M}_{R,\eta, t}\cap F$ } \hspace{2cm}
{\scriptsize $2\eta$} \hspace{8cm} {\scriptsize ${\cal M}_{R,\eta, t}\cap F$ }

\vspace{0,5cm}

{\scriptsize $M\cap F$}

\vspace{-0,4cm}

\hspace{2,7cm} {\scriptsize $y_0$}

\vspace{0,76cm}

\hspace{6,3cm} {\scriptsize $t_0$}

\vspace{-0,3cm}

\hspace{2,7cm} {\scriptsize $2\eta$}

\vspace{0,23cm}

\hspace{4,1cm} {\scriptsize $R$} \hspace{2,8cm} {\scriptsize $e$}

\vspace{0,2cm}

\hspace{6,1cm} {\scriptsize $t$}

\vspace{0,05cm}

\hspace{8,4cm} {\scriptsize $\xi$}

\vspace{0,1cm}

\hspace{7cm} {\scriptsize $O$}

{\scriptsize  $F$ is the vector space spanned by $e$ and $\xi$}


\end{proofcontactpoint}

\saut

\saut

\begin{proofnirmos} Integrating by parts we have

\begin{align*}\insm\frac{1}{2}\Delta\xi^2\xi^{2k-2}\vol&=\frac{1}{2}\insm\scal{d\xi^2}{d\xi^{2k-2}}\vol=2\left(\frac{k-1}{k^2}\right)\insm|d\xi^k|^2\vol\\
&\leq (A_1+kA_2)\insm\xi^{2k-r}\vol+(B_1+kB_2)\insm\xi^{2k-s}\vol\end{align*}

\ndt Now, given a smooth function $f$ and applying the Sobolev inequality
(\ref{sobol}) to $f^2$, we get

\begin{align*}\left(\insm f^{\frac{2n}{n-1}}\vol\right)^{1-(1/n)}&\leq K(n)\insm\left(2|f||df|+|H|f^2\right)\vol\\
&\leq 2K(n)\left(\insm f^2\vol\right)^{1/2}\left(\insm|df|^2\vol\right)^{1/2}+K(n)\nh\insm f^2\vol\displaybreak[2]\\
&=K(n)\left(\insm f^2\vol\right)^{1/2}\left(2\left(\insm|df|^2\vol\right)^{1/2}+\nh\left(\insm f^2\vol\right)^{1/2}\right)\end{align*}

\ndt where in the second inequality, we have used the H\"older
inequality. Using it again, by assuming that $V(M)=1$, we have 

$$\left(\insm f^2\vol\right)^{1/2}\leq\left(\insm f^{\frac{2n}{n-1}}\vol\right)^{\frac{n-1}{2n}}$$

\ndt And finally, we obtain

$$\|f\|_{\frac{2n}{n-1}}\leq K(n)\left(2\|df\|_2+\nh\|f\|_2\right)$$

\ndt For $k\geq 2$, $\xi^k$ is smooth and we apply the above inequality
to $f=\xi^k$. Then we  get

\begin{align*}\|\xi\|_{\frac{2kn}{n-1}}^k&\leq K(n)\left[2\left(\insm|d\xi^k|^2\vol\right)^{1/2}+\nh\left(\insm\xi^{2k}\vol\right)^{1/2}\right]\\
&\leq
K(n)\left[2\left(\frac{k^2}{2(k-1)}\right)^{1/2}\left((A_1+kA_2)\insm\xi^{2k-r}\vol+(B_1+kB_2)\insm\xi^{2k-s}\vol\right)^{1/2}\right.\displaybreak[2]\\
&\left.\hspace{9cm}+\nh\left(\insm\xi^{2k}\vol\right)^{1/2}\right]\\
&\leq K(n)\left[2\left(\frac{k^2}{2(k-1)}\right)^{1/2}
   \left((A_1+kA_2)\nxii^{2-r}+(B_1+kB_2)\nxii^{2-s}\right)^{1/2}\nfidkmd^{k-1}\right.
\\
&\hspace{9cm}\Biggl.+\nh\nxii\nfidkmd^{k-1}\ \ \Biggr]\displaybreak[2]\\
&\leq
K(n)\left[2\left(\frac{k^2}{2(k-1)}\right)^{1/2}\left(\frac{A_1+kA_2}{\nxii^r}+\frac{B_1+kB_2}{\nxii^s}\right)^{1/2}\right.
\Biggl.+\nh\Biggr]\nxii\nfidkmd^{k-1}\displaybreak[2]\\
&\leq K(n)\left[2\left(\frac{k^2}{2(k-1)}\right)^{1/2}
\left(\frac{A_1^{1/2}+k^{1/2}A_2^{1/2}}{\nxii^{r/2}}+\frac{B_1^{1/2}+k^{1/2}B_2^{1/2}}{\nxii^{s/2}}\right)\right.\\
&\hspace{9cm}\Biggl.+\nh\Biggr]\nxii\nfidkmd^{k-1}
\end{align*}

\ndt Now if we assume that $\nxii>\eta$, the last inequality becomes

\begin{align*}\|\xi\|_{\frac{2kn}{n-1}}^{k}&\leq
K(n)\left[2\left(\frac{k^2}{2(k-1)}\right)^{1/2}\left(\frac{A_1^{1/2}+k^{1/2}A_2^{1/2}}{\eta^{r/2}}+\frac{B_1^{1/2}+k^{1/2}B_2^{1/2}}{\eta^{s/2}}\right)\right.\\
&\Biggl.\hspace{9cm}+\nh\Biggr]\nxii\nfidkmd^{k-1}\displaybreak[2]\\
&=\left[(K_1+k^{1/2}K_2)\left(\frac{k^2}{k-1}\right)^{1/2}+K'\right]\nxii\nfidkmd^{k-1}\end{align*}

\ndt Now let $q=\frac{n}{n-1}>1$ and for $i\geq 0$ let $k=q^i+1\geq 2$ . Then

\begin{align*}\|\xi\|_{2(q^{i+1}+q)}&\leq\left(\left(K_1+(q^i+1)^{1/2}K_2\right)\left(\frac{q^i+1}{q^{i/2}}\right)+K''\right)^{\coefi}\nxii^{\coefi}\nfiqi\\
&\leq\left(\ktld q^{i}\right)^{\coefi}\nxii^{\coefi}\nfiqi\end{align*}

\ndt where $\ktld=2K_1+2^{3/2}K_2+K'$. We see that $\ktld$ has a finite limit when $\eta\rightarrow\infty$  and if $B_1 >0$, $\ktld\rightarrow\infty$ when
$\nh\rightarrow\infty$ or $\eta\rightarrow 0$. Moreover the H\"older inequality gives

$$\nfiqipi\leq\|\xi\|_{2(q^{i+1}+q)}$$

\ndt which implies

$$\nfiqipi\leq\left(\ktld q^{i}\right)^{\coefi}\nxii^{\coefi}\nfiqi$$

\ndt Now, by iterating from $0$ to $i$, we get

\begin{align*}\nfiqipi&\leq \ktld^{\left(1-\prod_{k=i-j}^i\left(1-\coefk\right)\right)}q^{\sum_{k=i-j}^i\frac{k}{q^k+1}}\nxii^{\left(1-\prod_{k=i-j}^i\left(1-\coefk\right)\right)}\|\xi\|_{2q^{i-j}}^{\prod_{k=i-j}^i\left(1-\coefk\right)}\\
&\leq \ktld^{\left(1-\prod_{k=0}^i\left(1-\coefk\right)\right)} q^{\sum_{k=0}^i\frac{k}{q^k+1}}\nxii^{\left(1-\prod_{k=0}^i\left(1-\coefk\right)\right)}\|\xi\|_{2}^{\prod_{k=0}^i\left(1-\coefk\right)}\end{align*}

\ndt Let $\alpha=\sum_{k=0}^{\infty}\frac{k}{q^k+1}$ and $\beta=\prod_{k=0}^{\infty}\left(1-\coefk\right)=\prod_{k=0}^{\infty}\left(\frac{1}{1+(1/q)^k}\right)$. Then

$$\nxii\leq
\ktld^{1-\beta}q^{\alpha}\nxii^{(1-\beta)}\nxid^{\beta}$$

\ndt and finally

$$\nxii\leq L\nxid$$

\ndt where $L=\ktld^{\frac{1-\beta}{\beta}}q^{\alpha/\beta}$ is a constant
depending only on $n$, $A_1$, $A_2$, $B_1$, $B_2$, $\nh$ and $\eta$. From classical methods we show
that $\beta\in [e^{-n},e^{-n/2}]$. In particular, $0<\beta<1$ and we deduce
that $L$ is bounded when $\eta\rightarrow\infty$ and $L\rightarrow\infty$ when $\nh\rightarrow\infty$ or $\eta\rightarrow
0$ with $B_1>0$.

\end{proofnirmos}

\saut

\ndt {\bf Remark} In \cite{shixu1} and \cite{shixu2} Shihohama and Xu have
proved that if $\var$ is  a compact $n$-dimensional Riemannian ma\-nifold
without boundary isometrically immersed in $\R^{n+1}$ and if
$\insm\left(|B|^2-n|H|^2\right)<D_n$ where $D_n$ is a constant depending on
$n$, then all Betti numbers are zero. For $n=2$, $D_2=4\pi$, and it
follows that if

$$\insm|B|^2\vol-4\pi<\la V(M)$$

\ndt then we deduce from the Reilly inequality $\la V(M)\leq 2\insm
H^2\vol$ that $\insm(|B|^2-2|H|^2)\vol<4\pi$ and by the result of Shihohama
and Xu $M$ is diffeomorphic to $\S^2$.

\end{document}